\documentclass[10pt,b5paper]{article}

\usepackage[
textwidth=125mm,
textheight=185mm,
centering
]{geometry}

\usepackage[T1]{fontenc}
\usepackage{lmodern}
\usepackage{amsmath,amssymb,amsthm,mathtools}
\usepackage{tabularx,microtype}

\usepackage{graphicx}
\usepackage{booktabs}
\usepackage{enumitem}
\usepackage[hidelinks]{hyperref}

\newtheorem{theorem}{Theorem}[section]
\newtheorem{lemma}[theorem]{Lemma}
\newtheorem{proposition}[theorem]{Proposition}
\newtheorem{corollary}[theorem]{Corollary}

\theoremstyle{definition}

\newtheorem{example}[theorem]{Example}

\theoremstyle{remark}
\newtheorem{remark}[theorem]{Remark}

\newcommand{\F}{\mathbb F}
\newcommand{\C}{\mathbb C}
\newcommand{\M}{\mathrm M}
\newcommand{\GL}{\mathrm{GL}}
\newcommand{\rk}{\operatorname{rank}}
\newcommand{\urk}{\operatorname{urk}}
\newcommand{\im}{\operatorname{im}}
\newcommand{\codim}{\operatorname{codim}}
\newcommand{\Hom}{\operatorname{Hom}}
\newcommand{\Gr}{\operatorname{Gr}}
\newcommand{\Alt}{\operatorname{Alt}}
\newcommand{\Span}{\operatorname{span}}
\newcommand{\D}{\mathcal D}
\newcommand{\Rcal}{\mathcal R}
\newcommand{\Pcal}{\mathcal P}
\newcommand{\T}{^{\mathsf T}}
\newcommand{\qbinom}[2]{\genfrac{[}{]}{0pt}{}{#1}{#2}_{Q}}
\newcommand{\doiurl}[1]{\href{https://doi.org/#1}{doi:\nolinkurl{#1}}}
\setlist[enumerate]{leftmargin=*,itemsep=2pt,topsep=4pt}
\title{Linear rank-\(k\) nonincreasing and preserving maps}

\author{Huajun Huang\thanks{Department of Mathematics and Statistics, Auburn University, Auburn, AL 36849, USA. Email: \texttt{huanghu@auburn.edu}}}

\date{}
\makeatletter
\renewcommand{\@maketitle}{%
  \begin{center}
    {\LARGE\bfseries\@title\par}
    \ifx\@author\@empty\else
      \medskip{\large\begin{tabular}[t]{c}\@author\end{tabular}\par}
    \fi
  \end{center}}
\makeatother

\begin{document}
	
\pagestyle{myheadings}
\markright{Rank-\(k\) nonincreasing and preserving maps}
\maketitle

\begin{abstract}
We characterize the linear maps between rectangular matrix spaces that carry matrices of rank at most \(k\) to matrices of rank at most \(k\). Over every field, such a map has the form \(X\mapsto MXN\) or \(X\mapsto MX\T N\), or its entire image has upper rank at most \(k\).
The proof uses polynomial functions on determinantal loci and extends Loewy's infinite-field theorem to all finite fields.
Over fields with more than two elements, exact preservation of rank \(k\) implies rank-\(k\) nonincrease. Over \(\F_2\), for each \(k\ge2\), we construct a rank-\(k\) preserver that sends a rank-one matrix to a matrix of rank \(k+1\).
For maps whose images have upper rank \(k\), three block rank equalities characterize exact preservation.
Over \(\C\), we combine the Huang--Landsberg classification with a degree obstruction on \(\mathrm{PGL}_k\) to show that every primitive core of upper rank four arising from a rank-four preserver is, up to equivalence and transposition, a projected exterior-product space on \(\C^5\).
If either source dimension exceeds five, every complex rank-four preserver is standard with full-rank factors or is a compression. A \(4\times5\) example attains the bound.
The degree obstruction also shows that, over any algebraically closed field, every linear map \(\M_k\to\Alt_{k+1}\) sends some invertible matrix to a matrix of rank at most \(k-2\) when \(4\mid k\).
\end{abstract}


\section{Introduction}\label{sec:main}
Let \(\phi:\M_{m,n}(\F)\to\M_{p,q}(\F)\) be a linear map between matrix spaces over a field \(\F\), and fix \(1\le k\le\min(m,n)\). We ask what the values of \(\phi\) on matrices of rank \(k\), or of rank at most \(k\), determine about \(\phi\) on the full matrix space. Write
\[
 \D_k^{m,n}=\{X:\rk X\le k\},\qquad
 \Rcal_k^{m,n}=\{X:\rk X=k\}
\]
for the determinantal locus and its rank-\(k\) stratum, respectively. We call \(\phi\) \emph{rank-\(k\) nonincreasing} if \(\phi(\D_k^{m,n})\subseteq\D_k^{p,q}\).
It \emph{preserves rank \(k\)} if \(\phi(\Rcal_k^{m,n})\subseteq\Rcal_k^{p,q}\), a condition requiring \(k\le\min(p,q)\). We suppress \(\F\) from the matrix-space notation and write \(\M_n=\M_{n,n}\). For a matrix space \(S\), its \emph{upper rank} is \(\displaystyle \urk S=\max_{Y\in S}\rk Y\), taken on matrices over the given field.

A \emph{standard map} has one of the forms
\begin{equation}\label{eq:standard}
  X\longmapsto MXN \qquad\text{or}\qquad X\longmapsto MX\T N
\end{equation}
with fixed factors of compatible sizes. In either case, \(\rk\phi(X)\le\rk X\) for every \(X\).

\begin{theorem}[Rank-\(k\) nonincrease over arbitrary fields]\label{thm:main}
The map \(\phi\) is rank-\(k\) nonincreasing if and only if it is standard or \(\urk(\im\phi)\le k\). Consequently every rank-\(k\) nonincreasing map satisfies
\[
  \rk\phi(X)\le\max\{k,\rk X\}\qquad\text{for every }X\in\M_{m,n}.
\]
\end{theorem}

Loewy proved the standard form for endomorphisms of a full matrix space over an algebraically closed field when the image contains a matrix of rank greater than \(k\) \cite{Loewy}. He later observed that the same conclusion holds over every infinite field \cite[Remark~5]{LoewyII}. For rank-\(k\) nonincrease, the rectangular version over an infinite field follows by taking a coordinate block of a sufficiently large square matrix, applying \(\phi\), and padding the output by zeros. Applying Loewy's result to this endomorphism and restricting back to the source block gives \eqref{eq:standard} whenever the image has upper rank greater than \(k\). When \(|\F|>k+1\), one can also pass to the algebraic closure. Write a matrix of rank at most \(k\) as \(X=AB\) with \(A\in\M_{m,k}\) and \(B\in\M_{k,n}\). Each \((k+1)\)-minor of \(\phi(AB)\) has degree at most \(k+1\) in each entry of \(A,B\), so vanishing over \(\F\) implies vanishing over its algebraic closure by interpolation \cite{HT}. If \(\urk(\im\phi)>k\), Loewy's theorem over the algebraic closure then gives rank-one nonincrease on \(\F\), and Botta's theorem \cite{Botta} gives one of the standard forms in \eqref{eq:standard} over \(\F\). Huang and Tsai classified rank-two nonincrease over fields with at least five elements using the Chow rank of the generic cubic determinant \cite{HT}.

The proof of Theorem~\ref{thm:main} works over the ground field, including \(\F_2\). Every polynomial of degree at most \(k+1\) that vanishes on \(\D_k^{m,n}\) agrees, as a function on \(\M_{m,n}(\F)\), with a linear combination of \((k+1)\)-minors. We apply this fact to square output blocks to show that \(\urk(\im\phi)>k\) forces rank-one nonincrease. Botta's theorem then gives the standard form.

In the bounded-image case, which spaces of upper rank \(k\) admit a linear parametrization that attains rank \(k\) on \emph{every} rank-\(k\) input? The Atkinson--Lloyd reduction \cite{AL} separates rows and columns whose deletion lowers the upper rank, leaving a primitive core after unused rows and columns are removed. When the core is zero, suitable output bases place the support of every image matrix in the union of \(a\) rows and \(b\) columns with \(a+b\le k\). We call a map with such support a \emph{compression} of upper rank at most \(k\). Over \(\C\), Huang and Landsberg classified the four basic spaces of upper rank four \cite{HL}. Requiring the parametrization to preserve rank four singles out one of their models.

\begin{theorem}[Rank-four preservers over \(\C\)]\label{thm:fourintro}
Suppose \(\F=\C\), \(k=4\), and \(\phi\) preserves rank four. If \(\max(m,n)>5\), then \(\phi\) is standard with a full-column-rank left factor and a full-row-rank right factor, or it is a compression with \(a+b=4\). The bound five is sharp.

For a bounded-image preserver, every primitive core of upper rank four is, up to equivalence and transposition, a projected exterior-product space on \(U=\C^5\). Its matrices have the form \(y\mapsto T(z\wedge y)\), where \(z\) ranges over \(U\) and the kernel of the linear projection \(T\) contains no nonzero decomposable bivector.
\end{theorem}

Theorem~\ref{thm:four} gives the full description of the possible cores. The exterior-product construction in Section~\ref{sec:ext}
attains the bound five on a \(4\times5\) source. To exclude the other basic models, we use linear conditions on their coefficients and the divisibility condition for the degree of a projective rank-drop locus proved in Proposition~\ref{prop:degree}. Corollary~\ref{cor:pfaffian} shows that, for every positive \(k\) divisible by four, each linear map \(\M_k\to\Alt_{k+1}\) over an algebraically closed field takes some invertible input to a matrix of rank at most \(k-2\).

Over a finite field \(\F_Q\), the passage from exact preservation to
rank nonincrease rests on Lemma~\ref{lem:finitevanish}. For \(k\ge2\),
a polynomial of degree less than \(k(Q-2)+2\) that vanishes on the
rank-\(k\) stratum also vanishes on \(\D_k^{m,n}\). When \(Q>2\), the
bound applies to every \((k+1)\)-minor of \(\phi(X)\).
Together with the matrix pencil of Lemma~\ref{lem:pencil}, this gives
rank-\(k\) nonincrease over every field with more than two elements
whenever \(\rk\phi(X)\le k\) on the rank-\(k\) stratum
(Corollary~\ref{cor:stratum}). Example~\ref{ex:binary} shows that the
field hypothesis is sharp for each \(k\ge2\).
Theorem~\ref{thm:pres} and Proposition~\ref{prop:blocks} characterize
exact preservation in the bounded-image case by three block rank
equalities. Section~\ref{sec:ext} describes the primitive space and its parametrization explicitly for exterior products.

Related exact-preserver results include Beasley's theorem for endomorphisms of a complex full matrix space \cite{Beasley}, the theorem of Beasley and Laffey for invertible endomorphisms over fields with at least four elements \cite{BL}, and the constructions of Li, Rodman and \v Semrl for maps between different rank strata \cite{LRS}. Huang and Tsai proved that exact rank-\(k\) preservation implies
rank-\(k\) nonincrease when \(|\F|>k+2\), and that rank-\(k\)
nonincrease implies rank-\(\ell\) nonincrease for each admissible
\(\ell\ge k\) when \(|\F|>k+1\) \cite{HT}. The global inequality in Theorem~\ref{thm:main} gives the latter implication without a field-size hypothesis.

\section{Polynomial functions on determinantal loci}\label{sec:minors}
We use the following lemma in place of interpolation. Its conclusion concerns the values of a polynomial on \(\F\)-points.

\begin{lemma}\label{lem:function}
Let \(r\ge2\). Suppose that a polynomial \(f(X)\) of total degree at most \(r\) vanishes on \(\D_{r-1}^{m,n}\). As a function on \(\M_{m,n}(\F)\),
\begin{equation}\label{eq:minor}
       f(X)=\sum_{\substack{|I|=r\\ |J|=r}}c_{I,J}\det X[I,J],
       \qquad c_{I,J}\in\F .
\end{equation}
The sum is empty when \(r>\min(m,n)\), and the coefficients are unique.
\end{lemma}
\begin{proof}
Group monomials according to the exact set \(J\) of columns appearing in them. Write \(X_H\) for \(X\) with columns outside \(H\) set to zero. The component with support \(J\) takes the value
\[
       f_J(X)=\sum_{H\subseteq J}(-1)^{|J|-|H|}f(X_H).
\]
For \(|J|<r\), every summand vanishes. For \(|J|=r\), the degree bound makes \(f_J\) multilinear in these \(r\) columns. If they are dependent, every matrix \(X_H\) in the sum has rank less than \(r\). Hence \(f_J\) is alternating and has a unique expansion in the coordinate alternating \(r\)-forms. These are precisely the minors with column set \(J\). No larger support occurs. To read off \(c_{I,J}\), insert an identity matrix at \(I,J\) and put zeros elsewhere.
\end{proof}

Over an infinite field, \eqref{eq:minor} is a polynomial identity and describes the degree-\(r\) part of the generic determinantal ideal \cite{BV}. Over a finite field, a nonzero polynomial can induce the zero function. The proof only uses values on field-valued matrices.

\begin{lemma}\label{lem:concise}
Over every field the determinant function on \(\M_r\) has no nonzero translation period. Any nonzero function of the form \eqref{eq:minor} requires at least \(r^2\) homogeneous linear coordinates, even if an arbitrary function of those coordinates is allowed.
\end{lemma}
\begin{proof}
For \(H\ne0\), change bases so that \(H=\operatorname{diag}(I_h,0)\). Taking \(Y=\operatorname{diag}(0_h,I_{r-h})\) gives \(\det Y=0\) and \(\det(Y+H)\ne0\). Invertible basis changes preserve this distinction. For the second assertion restrict to a single block \(I,J\) with \(c_{I,J}\ne0\). On that block the function is a nonzero multiple of the determinant. Fewer than \(r^2\) linear coordinates would have a nonzero common kernel on \(\M_r\), giving a translation period.
\end{proof}

A linear \(r\times r\) matrix with a fixed zero entry uses at most \(r^2-1\) linear coordinates. Its determinant therefore vanishes identically as a function if it vanishes on all inputs of rank less than \(r\).

\begin{proposition}\label{prop:squaretest}
Let a linear map \(\theta:\M_{m,n}\to\M_r\), \(r\ge2\), send \(\D_{r-1}^{m,n}\) into the singular matrices. If its image contains an invertible matrix, then \(\theta\) is a standard map. Both of its outer factors have rank \(r\).
\end{proposition}
\begin{proof}
Apply Lemma~\ref{lem:function} to \(\det\theta(X)\). Some coefficient in \eqref{eq:minor} is nonzero. After permuting source rows and columns, the restriction \(\theta_0:\M_r\to\M_r\) to the upper left block satisfies \(\det\theta_0(Y)=c\det Y\) for all \(Y\), with \(c\ne0\). Every element of \(\ker\theta_0\) is therefore a translation period of the determinant. By Lemma~\ref{lem:concise} this kernel is zero, so \(\theta_0\) is bijective. The invertible linear preserver theorem for \(\GL_r\) over arbitrary fields \cite[Theorem~1]{Inv} now gives \(\theta_0(Y)=PYQ\) or \(PY\T Q\), with \(P,Q\) invertible. Composing on the output with \(\theta_0^{-1}\), we may assume \(\theta_0(Y)=Y\).

Write the input as \(\left(\begin{smallmatrix}Y&B\\ C&D\end{smallmatrix}\right)\). For one additional column \(b\), write \(L(b)\) for its contribution. If \(H\) is a hyperplane containing \(b\), then \(Y+L(b)\) is singular for every \(Y\) with \(\im Y\subseteq H\). This forces \(\im L(b)\subseteq H\). Indeed, write \(H=\ker\lambda\). If \(\lambda L(b)\ne0\), one can extend that nonzero row to an invertible matrix \(Z\) satisfying \(\lambda Z=\lambda L(b)\), and take \(Y=Z-L(b)\). Intersecting over \(H\) gives \(\im L(b)\subseteq\F b\). Each column of the linear map \(b\mapsto L(b)\) is therefore a scalar multiple of \(b\), with the scalar independent of \(b\) by linearity on two basis vectors and their sum. The transpose argument applies to every additional row. Consequently,
\[
  \theta(X)=Y+BV+UC+L_0(D)
        =(\,I\ U\,)X\binom{I}{V}+H_0(D)
\]
for fixed \(U,V\) and a linear \(H_0\). Invertible block triangular source changes replace the first term by the new upper left block and leave \(D\) fixed. In these coordinates \(\theta(X)=Y+H_0(D)\).

For a matrix unit \(E_{ij}\) of the remaining block and any \(Y\) of rank \(r-1\), complete \(Y\) to an input of the same rank with \(D=E_{ij}\). To do this, take
\(\binom{I}{S}Y(\,I\ T\,)\), choosing just row \(i\) of \(S\) and column \(j\) of \(T\) nonzero so that their pairing through \(Y\) is \(1\). Thus \(Y+H_0(E_{ij})\) is singular for every such \(Y\). If \(H_0(E_{ij})\) has positive rank \(h\), put it in the form \(\operatorname{diag}(I_h,0)\) and choose \(Y=\operatorname{diag}(J_h,I_{r-h})\) in the same bases, with \(J_h\) a single nilpotent Jordan block. Then \(Y\) has rank \(r-1\) and the sum is invertible. This proves \(H_0=0\). Undoing the transformations gives \eqref{eq:standard}. An invertible output forces both factors to have rank \(r\).
\end{proof}

\begin{lemma}\label{lem:simul}
If \(A,B\in\M_{p,q}\), \(\rk A\ge s\), \(\rk B\ge r\), and \(1\le s\le r\), there are \(L\in\M_{r,p}\) and \(R\in\M_{q,r}\) such that \(\rk LAR\ge s\) and \(LBR\) is invertible.
\end{lemma}
\begin{proof}
Choose \(s\) domain vectors whose images are independent under both maps. At each step the forbidden vectors form the union of two proper linear subspaces, which cannot exhaust a vector space over any field. Extend to \(r\) vectors with independent \(B\)-images and use them as the columns of \(R\). Repeat this construction for \((AR)\T\) and \((BR)\T\) to select the rows of \(L\).
\end{proof}
\begin{proof}[Proof of Theorem~\ref{thm:main}]
Only the case \(\urk(\im\phi)>k\) needs proof. Put \(r=k+1\) and choose \(Y\) with \(\rk\phi(Y)\ge r\). If some rank-one \(X\) has \(\rk\phi(X)\ge2\), apply Lemma~\ref{lem:simul} with \(A=\phi(X), B=\phi(Y), s=2\). The map \(Z\mapsto L\phi(Z)R\) satisfies Proposition~\ref{prop:squaretest} and hence sends \(X\) to a matrix of rank at most one, a contradiction. Thus \(\phi\) is rank-one nonincreasing. Botta's theorem \cite{Botta} gives a standard form or an image of upper rank at most one. Our choice of \(Y\) excludes the latter possibility. Conversely, a standard map satisfies \(\rk\phi(X)\le\rk X\), and in the bounded-image case \(\rk\phi(X)\le k\) for every \(X\). These inequalities also give the asserted bound on \(\rk\phi(X)\).
\end{proof}

\section{From exact preservation to rank nonincrease}\label{sec:exact}
\begin{lemma}\label{lem:pencil}
Suppose \(|\F|>k+1\). If \(\rk\phi(X)\le k\) for every \(X\) of rank exactly \(k\), then \(\phi\) is rank-\(k\) nonincreasing.
\end{lemma}
\begin{proof}
Take \(A\) of rank \(h<k\). Choose a nilpotent \(k\times k\) matrix \(J\) of rank \(h\), for instance \(J_{h+1}\oplus0\). Choose fixed source bases so that \(A\) is this \(J\) padded by zeros. In the same bases let \(B\) be the padded \(I_k\). Both \(B\) and \(A+tB\), for \(t\ne0\), have rank \(k\). Each \((k+1)\)-minor of \(\phi(A)+t\phi(B)\) has degree at most \(k\), since its coefficient of \(t^{k+1}\) is a minor of \(\phi(B)\) and vanishes. Each such minor vanishes at more than \(k\) nonzero values of \(t\), so it also vanishes at \(t=0\). Hence \(\rk\phi(A)\le k\). If the output has fewer than \(k+1\) rows or columns the bound is automatic.
\end{proof}

\begin{lemma}\label{lem:digitaverage}
Let \(Q,k\ge2\) be integers, \(N=Q^k-1\), and \(D=N/(Q-1)\).
For \(0\le n\le N\), let \(s(n)\) be the sum of its \(k\) base-\(Q\)
digits. For \(m\ge1\), write \([e]_m\) for the representative of \(e\bmod m\) in
\(\{1,\ldots,m\}\). For every integer \(e\),
\[
 \frac{1}{Q-1}\sum_{t=0}^{Q-2}s([e+tD]_N)
 \ \ge\ \frac{k(Q-2)}{2}+1 .
\]
\end{lemma}
\begin{proof}
Multiplication by \(Q\) modulo \(N\) cyclically permutes the \(k\)
digits, including the case \([e]_N=N\). Summing the \(k\) rotations gives
\[
 D\,s([e]_N)=\sum_{h=0}^{k-1}[Q^h e]_N.
\]
For fixed \(h\), put \(r_h=[Q^h e]_D\). Since \(QD\equiv D\pmod N\),
as \(t\) runs from \(0\) to \(Q-2\), the numbers \([Q^h(e+tD)]_N\)
run through \(r_h,r_h+D,\ldots,r_h+(Q-2)D\). Averaging the rotation
identity therefore gives
\[
 \frac{1}{Q-1}\sum_{t=0}^{Q-2}s([e+tD]_N)
 =\frac{k(Q-2)}{2}+\frac{1}{D}\sum_{h=0}^{k-1}r_h .
\]
The last sum is positive and congruent to \(e\sum_{h=0}^{k-1}Q^h=eD\)
modulo \(D\), so it is at least \(D\).
\end{proof}

\begin{lemma}\label{lem:finitevanish}
Let \(\F=\F_Q\) and \(2\le k\le\min(m,n)\). If a polynomial in the
entries of \(X\in\M_{m,n}(\F_Q)\) has total degree less than
\(k(Q-2)+2\) and vanishes on every matrix of rank \(k\), then it
vanishes on every matrix of rank at most \(k\).
\end{lemma}
\begin{proof}
After row and column changes, any matrix of rank at most \(k\) is
supported in a \(k\times k\) block. Restriction to this block preserves
the degree bound and gives a polynomial \(f\) that vanishes on
\(\GL_k(\F_Q)\). It suffices to show that \(f\) vanishes on all of
\(\M_k(\F_Q)\).

Set \(K=\F_{Q^k}\) and use \(N,D,s\) as in Lemma~\ref{lem:digitaverage}.
Fix an \(\F_Q\)-basis of \(K\) and let \(L_{a,b}\) be the matrix of
\[
  u\longmapsto au-bu^Q,\qquad a,b\in K.
\]
If exactly one of \(a,b\) is zero, then \(L_{a,b}\) is invertible. For
\(ab\ne0\), it is singular precisely when \(a/b=u^{Q-1}\) for some
\(u\in K^\times\), or equivalently \(a^D=b^D\), since \(K^\times\)
is cyclic. In particular, \(\ker L_{1,1}=\F_Q\), so \(L_{1,1}\) has
rank \(k-1\).

Represent the function \(g(a,b)=f(L_{a,b})\) by its unique polynomial
over \(K\) of degree at most \(N\) in each variable. Every monomial
\(a^i b^j\) with nonzero coefficient satisfies
\[
 s(i)+s(j)\le\deg f .
\]
Indeed, the entries of \(L_{a,b}\) are \(\F_Q\)-linear in \(a,b\).
By nondegeneracy of the trace pairing, every \(\F_Q\)-linear
functional on \(K\) has the form
\(a\mapsto\operatorname{Tr}_{K/\F_Q}(ca)=\sum_{h=0}^{k-1}(ca)^{Q^h}\).
In expanding a monomial of \(f\), each factor thus contributes one
power \(a^{Q^h}\) or \(b^{Q^h}\). To reduce the exponent of \(a\),
replace \(Q\) copies of \(a^{Q^h}\) by one copy of \(a^{Q^{h+1}}\),
using \(a^{Q^k}=a\) when \(h=k-1\). Each replacement decreases the
number of factors. Once every multiplicity is less than \(Q\), the
multiplicities are the digits of the reduced exponent. The same
argument applies to \(b\).

On the \(a\)-axis the reduced polynomial is \(g(0,0)(1-a^N)\).
Since \(s(N)=k(Q-1)\ge k(Q-2)+2\), the digit bound forces the coefficient of \(a^N\) to vanish.
It follows that \(g(0,0)=0\). Thus \(g\) vanishes on both axes, and
each monomial with nonzero coefficient has \(1\le i,j\le N\).

On \((K^\times)^2\), the support of \(g\) lies in \(a^D=b^D\), so
\(a^D b^{-D}g=g\). The monomials \(a^ib^j\), \(1\le i,j\le N\),
form a basis for the \(K\)-valued functions on \((K^\times)^2\).
Hence a nonzero coefficient at \((i,j)\) forces a nonzero coefficient
at each pair
\[
 \bigl([i+tD]_N,[j-tD]_N\bigr),\qquad 0\le t\le Q-2 .
\]
Since \(tD\) and \(-tD\) run through the same residues modulo \(N\),
Lemma~\ref{lem:digitaverage} bounds the average of
\(s([i+tD]_N)+s([j-tD]_N)\) from below by \(k(Q-2)+2\).
Each pair with nonzero coefficient has digit sum at most
\(\deg f<k(Q-2)+2\). Thus \(g=0\).

In particular \(f(L_{1,1})=0\). Applying the argument to
\(X\mapsto f(PXR)\) for \(P,R\in\GL_k(\F_Q)\) shows that \(f\)
vanishes on every matrix of rank \(k-1\). Now suppose that \(f\)
vanishes on every matrix of rank \(r\), with \(1\le r<k\).
For any \(A\) of rank \(r\), the polynomial \(X\mapsto f(AX)\) vanishes
on \(\GL_k(\F_Q)\), so \(f(AB)=0\) for every \(B\) of rank \(k-1\).
Every matrix of rank \(r-1\) is such a product, as follows from
\(\operatorname{diag}(I_r,0)(I_k-E_{rr})=\operatorname{diag}(I_{r-1},0)\)
by changes of bases. Descending induction completes the proof.
\end{proof}

\begin{corollary}\label{cor:stratum}
Suppose \(k=1\) or \(|\F|>2\). If \(\rk\phi(X)\le k\) whenever
\(\rk X=k\), then \(\phi\) is rank-\(k\) nonincreasing.
\end{corollary}
\begin{proof}
For \(k=1\) use \(\phi(0)=0\). Over an infinite field apply
Lemma~\ref{lem:pencil}. Over \(\F_Q\) with \(Q>2\) and \(k\ge2\), each
\((k+1)\)-minor of \(\phi(X)\) has degree at most \(k+1<k(Q-2)+2\)
and vanishes whenever \(\rk X=k\). By Lemma~\ref{lem:finitevanish}, each of these minors vanishes on \(\D_k^{m,n}\).
\end{proof}

\begin{theorem}[Exact rank-\(k\) preservation]\label{thm:pres}
Suppose \(k\le\min(p,q)\) and that either \(k=1\) or \(|\F|>2\). Then \(\phi\) preserves rank \(k\) if and only if at least one of the following holds.
\begin{enumerate}[label=\textup{(\roman*)}]
\item It is standard, its left factor has full column rank, and its right factor has full row rank.
\item Its image has upper rank \(k\). In the block form \eqref{eq:AL} below, \(a+b+r=k\), and the three block equalities \eqref{eq:three} hold on every rank-\(k\) input.
\end{enumerate}
The same characterization holds over any field for maps already known to be rank-\(k\) nonincreasing.
\end{theorem}
\begin{proof}
By Corollary~\ref{cor:stratum}, every rank-\(k\) preserver under the stated hypotheses is rank-\(k\) nonincreasing. Apply Theorem~\ref{thm:main}. Consider first the standard form \(X\mapsto MXN\). If \(\ker M\ne0\), some rank-\(k\) input has a column space meeting \(\ker M\) nontrivially, so its rank drops. The same argument on row spaces forces \(N\) to have full row rank. Since transposition preserves rank, both conclusions also hold for \(X\mapsto MX\T N\). Conversely, factors of these ranks preserve the rank of every input. In the bounded case the upper rank must be exactly \(k\). Proposition~\ref{prop:blocks} below gives the necessary and sufficient conditions for attaining it.
\end{proof}

\begin{example}\label{ex:binary} 
	For every \(k\ge2\), there is a linear map
	\(\phi:\M_k(\F_2)\to\M_{k+1}(\F_2)\) that preserves rank \(k\)
	and sends a rank-one matrix to a matrix of rank \(k+1\).
	
	Fix an \(\F_2\)-basis of \(K=\F_{2^k}\), and use it to identify
	\(\F_2^k\) with \(K\). For \(X=(x_{ij})\in\M_k(\F_2)\), let
	\(v\in K\) represent the first column of \(X\), and put \(t=x_{12}\).
	Define \(\phi(X)\) on \(\F_2\oplus K\) by
	\[
	\phi(X)(s,u)=(ts,\,vu+tu^2).
	\]
	Since squaring is \(\F_2\)-linear on \(K\), this defines an
	\(\F_2\)-linear operator, and its dependence on \(X\) is linear.
	
	If \(X\) is invertible, then \(v\ne0\). For \(t=0\), multiplication
	by \(v\) is invertible on \(K\), so \(\phi(X)\) has rank \(k\).
	For \(t=1\), the equation \(vu+u^2=u(u+v)=0\) has precisely the
	solutions \(0\) and \(v\). Thus the second component has rank \(k-1\),
	and the first has rank one. Hence \(\phi\) preserves rank \(k\).
	
	For the rank-one input \(X=E_{12}\), we have \(v=0\) and \(t=1\).
	Since squaring is an automorphism of \(K\), the map \(\phi(E_{12})(s,u)=(s,u^2)\) is invertible.
\end{example}

\section{Bounded images and their parametrizations}\label{sec:blocks}
\subsection{The primitive core}
A space \(S\subseteq\Hom(\F^q,\F^p)\) is \emph{reduced} if its common kernel is zero and its images span \(\F^p\). For \(r=\urk S>0\), a reduced space is \emph{primitive} if its upper rank stays \(r\) both on restriction to every domain hyperplane and after passage to every quotient of the codomain by a line. In coordinates, after any fixed change of the two output bases, deleting a row or a column still leaves some member of rank \(r\). The rank-\(r\) matrix may depend on which row or column is deleted. We allow the zero space of \(0\times0\) matrices as a trivial core. Applying the decomposition of Atkinson and Lloyd \cite[Theorem~1]{AL} to the image of \(\phi\) gives the following block form.

\begin{theorem}\label{thm:al}
Set \(s=\urk(\im\phi)\). There are invertible output factors \(U_0,V_0\) and nonnegative \(a,b,u,v,r\) such that
\begin{equation}\label{eq:AL}
 \Phi(X)=U_0\phi(X)V_0=
 \begin{pmatrix}A(X)&B(X)&E(X)\\ C(X)&P(X)&0\\ F(X)&0&0\end{pmatrix},
 \qquad a+b+r=s .
\end{equation}
The row sizes are \(a,u,p-a-u\), the column sizes \(b,v,q-b-v\). The displayed blocks are linear and \(P\) maps onto a primitive space \(\Pcal\subseteq\M_{u,v}\) of upper rank \(r\). Conversely, arbitrary such blocks give \(\rk\Phi(X)\le a+b+r\). When \(|\F|>r>0\), \(u,v\le\binom{r+1}{2}\).
\end{theorem}
\begin{proof}
The block reduction works over any field. After a change of bases, remove a common zero row or column whenever one exists. If the remaining space is reduced and imprimitive, delete a row or column whose removal lowers the upper rank. A single deletion lowers it by exactly one. Repeat until the remaining space is primitive or zero. Reversing the steps places the rank-dropping rows at the top and columns at the left, and the common zero rows and columns at the bottom and right. This produces \eqref{eq:AL} by fixed operations, and the block coordinates are linear maps. Conversely, removal of the first \(a\) rows and \(b\) columns leaves \(P(X)\) padded by zeros. When \(|\F|>r>0\), the bounds on \(u\) and \(v\) follow by applying \cite[Theorem~1]{AL} to the primitive space. The choice of core need not be unique.
\end{proof}

\begin{proposition}\label{prop:blocks}
For a fixed output in \eqref{eq:AL}, put
\[
 R=\begin{pmatrix}B&E\\P&0\end{pmatrix},\qquad
 L=\begin{pmatrix}C&P\\F&0\end{pmatrix}.
\]
If \(\rk P\le r\), then \(\rk\Phi=a+b+r\) if and only if
\begin{equation}\label{eq:three}
            \rk P=r,\qquad \rk R=r+a,\qquad \rk L=r+b .
\end{equation}
The block \(A\) is unrestricted. For \(s=k\), the equalities on all \(X\in\Rcal_k^{m,n}\) characterize exact preservation over any field.
\end{proposition}
\begin{proof}
The upper bound \(a+b+\rk P\) first forces \(\rk P=r\). Put \(P=\operatorname{diag}(I_r,0)\) in suitable bases. Removing its pivot rows and columns by elimination leaves
\[
 \begin{pmatrix}A-B_1C_1&\widehat B\\\widehat C&0\end{pmatrix},
 \qquad \widehat B=(B_0\ E),\qquad\widehat C=\binom{C_0}{F},
\]
where \(B_1,C_1\) meet the pivots of \(P\). Its rank is at most both \(b+\rk\widehat B\) and \(a+\rk\widehat C\). It attains \(a+b\) exactly when \(\widehat B\) has full row rank \(a\) and \(\widehat C\) has full column rank \(b\). In that case the columns of the right block span the first \(a\) coordinate directions, and the first \(b\) columns are independent modulo their span. Finally \(\rk R=r+\rk\widehat B\) and \(\rk L=r+\rk\widehat C\).
\end{proof}

For \(r=0\) the output is supported on \(a\) rows and \(b\) columns. Write \(\mathcal C_{a,b}\) for the full space with this support. A map into \(\mathcal C_{a,b}\) with \(a+b=k\) has the form
\begin{equation}\label{eq:compression}
    \begin{pmatrix}A(X)&B(X)\\ C(X)&0\end{pmatrix}.
\end{equation}
It preserves rank \(k\) exactly when \(B(X)\) has rank \(a\) and \(C(X)\) has rank \(b\) on that stratum. For a general bounded-image map, the condition \(\ker\phi\cap\Rcal_k=\varnothing\) alone does not ensure exact preservation. For example, a coordinate permutation of \(\M_k\) moving the diagonal to the first row is bijective and sends \(I_k\) to rank one for \(k\ge2\).

\subsection{Dimension bounds}
\begin{lemma}\label{lem:codim}
Put \(t=\max(m,n)\). If a linear map \(\Gamma:\M_{m,n}\to W\) satisfies \(\Gamma^{-1}(H)\cap\Rcal_k^{m,n}=\varnothing\) for a linear subspace \(H\le W\), then
\begin{equation}\label{eq:codim}
            t\le\dim\big((\im\Gamma+H)/H\big)\le\codim_W H .
\end{equation}
\end{lemma}
\begin{proof}
Restrict the quotient map to a \(k\times t\) source block, transposing the input if needed. Its kernel contains only matrices of rank at most \(k-1\). By Flanders's bound, the kernel has dimension at most \((k-1)t\), so the quotient map has rank at least \(t\). See \cite{Flanders} for the original theorem and \cite[Theorem~1.1]{Large} for the version over arbitrary fields.
\end{proof}

For a bounded-image preserver, every \(H\subseteq\Pcal\) of upper rank \(<r\) gives \(t\le\codim_{\Pcal}H\), by \eqref{eq:three}. More generally this holds in any linear parameter space in which membership in \(H\) forces a drop of core rank, whether or not the parameters are all used. In particular \(t\le\dim\Pcal\) if \(r>0\). If \(r>0\) and \(|\F|>r\), Theorem~\ref{thm:al} and Flanders's bound give
\begin{equation}\label{eq:large}
        t\le r\max(u,v)\le r\binom{r+1}{2}.
\end{equation}
Consequently if \(|\F|>k\) and \(t>k\binom{k+1}{2}\), a preserver is standard with full factors or admits \eqref{eq:compression}. For a compression the outer block sizes satisfy \(q-b\ge t\) if \(a>0\), and \(p-a\ge t\) if \(b>0\). Indeed, apply Lemma~\ref{lem:codim} to \(X\mapsto u\T B(X)\) or \(X\mapsto C(X)v\) with a nonzero fixed \(u\) or \(v\).

Flanders's theorem gives \(\dim S\le s\max(p,q)\) for \(S=\im\phi\) of upper rank \(s\). At equality, a linear space of positive upper rank \(s<\min(p,q)\) is a full one-sided compression, up to equivalence and transposition \cite[Theorem~1.1]{Large}. The following refinement is due to Atkinson and Lloyd in the square case, to Beasley in the rectangular case, and to de Seguins Pazzis over arbitrary fields \cite{ALLarge,Large}. Transpose the output if necessary so that \(p\ge q\). If \(1\le s<q\) and
\[
              \dim S\ge ps-(p-q+s)+1 ,
\]
then \(S\) is contained, up to equivalence, in one of
\(\mathcal C_{s,0},\mathcal C_{s-1,1},\mathcal C_{1,s-1},\mathcal C_{0,s}\), with the single possible exception \((p,q,s,\F)=(3,3,2,\F_2)\). The exceptional space is the five-dimensional space of upper triangular trace-zero matrices \cite[Theorem~1.3]{Large}. Thus in this dimension range \eqref{eq:compression} tests bounded-image preservation directly whenever the compression alternative holds.

\section{Projected exterior products}\label{sec:ext}
Let \(k\ge2,\ U=\F^{k+1}\) and \(N=\binom{k+1}{2}\). The classical exterior-product space \cite{Atkinson,HL} is given by the linear matrix \(W(v):U\to\bigwedge^2 U\), \(y\mapsto v\wedge y\). For \(v\ne0\), its kernel is \(\F v\) and its rank is \(k\). Fix a surjection \(T:\bigwedge^2U\to\F^d\) and a linear \(\Gamma:\M_{m,n}\to U\). Consider
\begin{equation}\label{eq:wedge}
              \phi_{T,\Gamma}(X)=T W(\Gamma(X))\in\M_{d,k+1}.
\end{equation}

\begin{theorem}\label{thm:wedge}
Over every field, \eqref{eq:wedge} preserves rank \(k\) if and only if \(\ker T\) contains no nonzero decomposable bivector and \(\ker\Gamma\) contains no matrix of rank \(k\).
\end{theorem}
\begin{proof}
The first condition makes \(T\) injective on \(w\wedge U\) for every \(w\ne0\), so the two conditions suffice. Conversely, if \(\Gamma(X)=0\) for a rank-\(k\) input, then \(\phi_{T,\Gamma}(X)=0\). If \(0\ne u\wedge v\in\ker T\), put \(H=\Span(u,v)\). For every nonzero \(w\in H\), the space \(w\wedge U\) contains \(u\wedge v\), and hence \(\rk TW(w)<k\). The rank also drops at \(w=0\). But \(H\) has codimension \(k-1\) in \(U\), and Lemma~\ref{lem:codim} gives a rank-\(k\) input in \(\Gamma^{-1}(H)\). This proves necessity even when \(\Gamma\) is not surjective.
\end{proof}

\begin{proposition}\label{prop:coeff}
Assume \(m,n\ge k\ge2\). A linear \(\Gamma:\M_{m,n}\to\F^{k+1}\) can avoid rank \(k\) in its kernel only for
\[
                  (m,n)=(k,k+1),\quad(k+1,k),\quad(k,k).
\]
Up to invertible changes of coefficient coordinates, the following list gives exactly the maps with this property.
\begin{enumerate}[label=\textup{(\alph*)}]
\item In the two rectangular cases, \(\Gamma(X)=X\T u\) or \(\Gamma(X)=Xv\), respectively, for a fixed nonzero vector of length \(k\).
\item On \(\M_k\), the possibilities include \(\Gamma(X)=(Xv,\ell(X))\) and \((X\T u,\ell(X))\), where \(\ell\) is any scalar linear form and \(u,v\ne0\). For \(k\ne3\) these exhaust the square case.
\item On \(\M_3\) there is also the quotient by \(\mathcal C_{1,1}\), given by the four entries of the complementary \(2\times2\) block. Over \(\F_2\) there is one further kernel type, the space of upper triangular trace-zero matrices. These descriptions allow fixed row and column changes on the source. In each quotient case any identification of the four-dimensional quotient with \(\F^4\) is allowed.
\end{enumerate}
\end{proposition}
\begin{proof}
Lemma~\ref{lem:codim} gives \(\max(m,n)\le k+1\). In \(\M_{k,k+1}\) the kernel has dimension at least \((k-1)(k+1)\) and upper rank at most \(k-1\). Equality in Flanders's bound forces it to be the space of maps into a fixed hyperplane of \(\F^k\) \cite[Theorem~1.1]{Large}. This proves (a) for \(k\times(k+1)\) matrices. Transposition gives the other case.

In size \((k+1)\times(k+1)\), write \(\Gamma(X)=\sum_i B_i x_i\) in terms of the transposed rows \(x_i\). Restricting to any \(k\) rows gives (a), so within that set the \(B_i\) are scalar multiples of one invertible matrix. Any two rows occur together because \(k\ge2\). Fix a nonzero \(B_i\). Then \(\Gamma(X)=B_i X\T c\) for a nonzero \(c\). Taking \(\im X=\ker c\T\) gives a rank-\(k\) input killed by \(\Gamma\).

In size \(k\times k\), the singular space \(K=\ker\Gamma\) has dimension at least \(k^2-k-1\). For \(k\ge3\) this is at least \(k(k-1)-(k-1)+1\), so the large-space theorem \cite[Theorem~1.3]{Large} applies. The two mixed compressions with parameters \(1,k-2\) and \(k-2,1\) have dimension \(k^2-2k+2\), too small if \(k>3\). At \(k=3\), equality forces the full mixed compression or the stated binary exception. Otherwise \(K\) has a common kernel vector or a common left kernel vector. The same conclusion holds for \(k=2\), because a space of rank at most one has a common image line or a common row line. Indeed the sum of two rank-one matrices with both pairs of factors independent has rank two. Hence the evaluation \(Xv\) or \(X\T u\) factors through \(\Gamma\), and its \(k\) independent coordinates can be completed to coordinates as in (b). Conversely, in the square case each listed kernel consists of singular matrices. For a full-rank rectangular input, the map \(X\T\) or \(X\) in (a) is injective, so the indicated evaluation is nonzero.
\end{proof}

In particular \(W(X\T u)\) on \(\M_{k,k+1}\) preserves rank \(k\) and takes some rank-one input to rank \(k\). This gives a nonstandard preserver in every rank \(k\ge2\), over every field.

\begin{proposition}\label{prop:wprimitive}
Over an algebraically closed field, suppose \eqref{eq:wedge} preserves rank \(k\ge2\). Its displayed matrix space is primitive if and only if \(\Gamma\) is surjective.
\end{proposition}
\begin{proof}
If \(\im\Gamma\) lies in a hyperplane \(H\subset U\), restriction of the matrices to \(H\) reduces the upper rank to at most \(k-1\), since \(z\wedge H\) has dimension at most \(k-1\) for \(z\in H\). Now suppose \(\Gamma\) is onto. The kernels \(\F z\) of the nonzero matrices intersect trivially, and their images span the target because decomposable bivectors span \(\bigwedge^2 U\). For any hyperplane \(H\subset U\), taking \(z\notin H\) gives rank \(k\) on restriction to \(H\).

A row deletion reducing the upper rank would give a line \(L\) in every \(T(z\wedge U)\), \(z\ne0\). In \(\mathbb P(T^{-1}(L))\), the decomposable points form a closed projective set disjoint from \(\mathbb P(\ker T)\) by Theorem~\ref{thm:wedge}. If \(\ker T\ne0\), then \(\mathbb P(\ker T)\) is a hyperplane in \(\mathbb P(T^{-1}(L))\). The set of decomposable points is therefore finite by the projective dimension theorem. If \(\ker T=0\), the ambient space is itself a point. Each decomposable point corresponds to a two-dimensional subspace of \(U\), and it belongs to \(\mathbb P(z\wedge U)\) only when that subspace contains \(z\). A finite union of such subspaces cannot cover \(U\), of dimension \(k+1\ge3\) over an infinite field. Thus no such line \(L\) exists.
\end{proof}

The minimum target dimension in the next proposition is the invariant
\(s_{k+1}(\F)\) of Gow and Quinlan. They determine its value over
finite fields and record the algebraically closed value in characteristic
different from two \cite[\S1, Theorem~2 and Corollary~8]{GQ}.

\begin{proposition}\label{prop:wtarget}
The least possible \(d\) for a surjective \(T\) as in Theorem~\ref{thm:wedge} is \(2k-1\) over an algebraically closed field and \(k+1\) over a finite field. Every \(d\) from that minimum to \(N\) occurs.
\end{proposition}
\begin{proof}
Over an algebraically closed field the decomposable lines form \(\Gr(2,k+1)\subseteq\mathbb P^{N-1}\), of dimension \(2k-2\). A nonempty projective center of codimension \(d\) must meet it if \(d\le2k-2\). For \(d>2k-2\) a general such center is disjoint, by the projective dimension theorem and the incidence correspondence \cite[Ch.~I]{Hartshorne}. The zero kernel is allowed.

Over \(\F_Q\), identify \(U\) with \(K=\F_{Q^{k+1}}\). The alternating map \((u,v)\mapsto uv^Q-u^Qv\) vanishes only on dependent pairs, because \(w^Q=w\) in \(K\) precisely for \(w\in\F_Q\). Taking the image of this map as the target gives \(d\le k+1\). Conversely \(d\ge k\), since \(T\) must be injective on \(w\wedge U\). If \(d=k\), every nonzero form in the \(k\)-dimensional image of \(T^*:(\F_Q^k)^*\to\bigwedge^2 U^*\) would be nondegenerate, since for \(w\ne0\) the map \(T(w\wedge-)\) is onto. Alternating forms in odd dimension are singular. In even dimension \(k+1=2e\) the Pfaffian on this space is a homogeneous polynomial of degree \(e<k\), so it has a nonzero zero over \(\F_Q\). Indeed, for a polynomial \(g\) of degree \(<k\) in \(k\) variables, summing \(1-g^{Q-1}\) over the field gives zero in \(\F_Q\) because every monomial has some exponent \(<Q-1\), and the sum over that variable is zero. Thus a homogeneous \(g\) cannot have precisely one zero. The Pfaffian argument also holds in characteristic two for alternating forms. Finally, shrinking any admissible kernel realizes all larger \(d\).
\end{proof}

\section{Primitive cores in low rank}\label{sec:low}
\subsection{Upper ranks one, two, and three}
A nonzero upper-rank-one space has a common image line or a common row line. It has no primitive core of positive rank. For \(|\F|>2\) the only primitive upper-rank-two space is \(\Alt_3(\F)\), the alternating \(3\times3\) matrices, including the zero-diagonal requirement in characteristic two \cite{AL,Pazzis}. Over \(\F_2\) all primitive upper-rank-two spaces still have format \(3\times3\). They are equivalent to a primitive subspace of the upper triangular trace-zero space, to \(\Alt_3\), or to one of the two explicit spaces \(U_3,V_3\) of de Seguins Pazzis \cite[Theorem~1.5]{Binary}.
For \(k=2\) and \(|\F|>2\), Theorem~\ref{thm:pres} gives the following description of rank-two preservers. Besides standard maps with full factors and compressions satisfying \eqref{eq:three}, the remaining maps have the full \(\Alt_3\) core with no additional row or column. In exterior-product coordinates they are \(X\mapsto W(\Gamma(X))\), up to equivalence, transposition and zero padding. Here \(\Gamma\) is surjective and its kernel contains no rank-two matrix. Proposition~\ref{prop:coeff} lists the possible maps \(\Gamma\). Indeed every nonzero bivector in \(\bigwedge^2\F^3\) is decomposable, so the admissible \(T\) in Theorem~\ref{thm:wedge} is an isomorphism.

For \(|\F|>3\), Theorem~C of Atkinson \cite{Atkinson}, in the arbitrary-characteristic form of \cite[Theorem~5.8]{Pazzis}, says that every primitive upper-rank-three space is, up to equivalence and transposition, obtained from
\[
       z\longmapsto T(z\wedge-),\qquad
       T:\bigwedge^2\F^4\twoheadrightarrow\F^d,\quad d=4,5,6 .
\]
All projections give bounded images. Primitivity imposes further conditions on the projection. Thus for \(k=3,\ |\F|>3\), Theorem~\ref{thm:main} gives standard maps, compressions with \(a+b\le3\), an alternating \(3\times3\) core with at most one extra row or column, and projected exterior products with no extra row or column. Zero padding and fixed output equivalences are understood. Every displayed family is rank-three nonincreasing.

For exact rank three, the compression test is \eqref{eq:compression}. With a rank-two core one needs \(a+b=1\). In the extra-row case the nonzero part can be written
\[
  \begin{pmatrix}h(X)&\ell(X)\\ C(z(X))&0\end{pmatrix},
  \qquad C(z)=\begin{pmatrix}0&-z_3&z_2\\z_3&0&-z_1\\-z_2&z_1&0\end{pmatrix}.
\]
For \(z\ne0\), \(C(z)\) has rank two and kernel \(\F z\). The output has rank three exactly when \(z\ne0\) and either \(\ell\ne0\) or \(hz\ne0\). Transposing gives the condition for an extra column. Lemma~\ref{lem:codim} forces \(m=n=3\), and equality in Flanders's theorem makes \(z\) an invertible change of \(Xv\) or \(X\T u\). With a rank-three core, Theorem~\ref{thm:wedge} and Proposition~\ref{prop:coeff} describe the projection and all coefficient maps. In Pl\"ucker coordinates a nonzero bivector in \(\bigwedge^2\F^4\) is decomposable precisely when
\[
            w_{12}w_{34}-w_{13}w_{24}+w_{14}w_{23}=0 .
\]
Hence \(d=4\) requires this quadratic form to have no nonzero zero on the two-dimensional kernel. This occurs over finite fields, by Proposition~\ref{prop:wtarget}, and cannot occur over \(\C\). 
Over every field with more than three elements, the preceding bounded-image families and the standard maps with full-rank factors exhaust the rank-three preservers. The hypothesis \(|\F|>3\) enters through the classification
of primitive spaces. Corollary~\ref{cor:stratum} also gives rank-three
nonincrease over \(\F_3\).

\subsection{Degrees of rank-drop loci}
Restricting the source of a rank-\(k\) preserver to a coordinate \(k\times k\) block gives a map that attains rank \(k\) on every invertible input. The following divisibility condition tests whether the coefficient map of a proposed model can avoid its rank-drop locus on this block.

\begin{proposition}\label{prop:degree}
Let \(\F\) be algebraically closed, \(1\le j<k\), and let \(Z\subset\mathbb P(W)\) be a projective subvariety of codimension \(j\) and degree \(\delta\). Suppose a linear map \(\Gamma:\M_k\to W\) avoids the affine cone over \(Z\), including its vertex, on \(\GL_k\). Then
\[
             \gcd\!\left\{\binom{k}{1},\binom{k}{2},\ldots,\binom{k}{j}\right\}\ \mid\ \delta .
\]
\end{proposition}
\begin{proof}
Let \(B=\mathbb P(\M_k)\setminus\{\det=0\}=\mathrm{PGL}_k\), and let \(h=c_1(\mathcal O_B(1))\). We compute the order of \(h^j\) in the integral Chow group \(\mathrm{CH}^j(B)\), the group of codimension-\(j\) cycles modulo rational equivalence. For the Chow-ring computation over \(\C\), see \cite[\S3]{DL}. We compute the order over the given field by localization and the projective bundle formula \cite{Fulton}.

In \(\mathbb P(\M_k)\times\mathbb P^{k-1}\) take
\[
       Y=\{([X],[v]):Xv=0\}.
\]
It is a projective bundle over \(\mathbb P^{k-1}\), with fibers \(\mathbb P^{k(k-1)-1}\). If \(H,x\) are the two hyperplane classes, its class in the product is \((H+x)^k\), since it is cut out in the expected codimension by \(k\) equations of bidegree \((1,1)\). Every integral subvariety of the determinant hypersurface lifts birationally to an integral subvariety of \(Y\). Indeed, over the function field of the subvariety, its generic matrix has a nonzero kernel vector. The line spanned by this vector defines a rational section, whose image closure gives the required lift. Consequently the cycles supported on the determinant hypersurface have the same image in the Chow group of \(\mathbb P(\M_k)\) as proper pushforwards from \(Y\).

In codimension \(j-1\) on \(Y\) the projective bundle formula gives generators \(H^{j-1-i}x^i\), \(0\le i<j\). To push each generator to \(\mathbb P(\M_k)\), multiply by \((H+x)^k\) and take the coefficient of \(x^{k-1}\). The result is \(\binom{k}{i+1}H^j\). Localization therefore identifies \(\mathrm{CH}^j(B)\) with \(\mathbb Z/b_{k,j}\mathbb Z\), generated by \(h^j\), where \(b_{k,j}=\gcd_{1\le i\le j}\binom{k}{i}\).

The linear system defining \(\Gamma\) has no base point on \(B\). It therefore defines a morphism \(B\to\mathbb P(W)\setminus Z\) under which the hyperplane class pulls back to \(h\). In the Chow ring of this open target, \(\delta H^j=0\), since \([Z]=\delta H^j\) in its projective completion. Pullback between the smooth open varieties gives \(\delta h^j=0\), proving the divisibility.
\end{proof}

\begin{corollary}[Pfaffian obstruction]\label{cor:pfaffian}
Let \(\F\) be algebraically closed. Write \(\Alt_n(\F)\) for the space of alternating \(n\times n\) matrices, with zero diagonal in every characteristic. If \(k\ge4\) is divisible by four, every linear map
\[
              \Gamma:\M_k(\F)\longrightarrow\Alt_{k+1}(\F)
\]
takes some invertible matrix to a matrix of rank at most \(k-2\).
\end{corollary}
\begin{proof}
Put \(k=2d\). The projective locus \(Z\) of alternating matrices of size \(2d+1\) and rank at most \(2d-2\) is defined by the submaximal Pfaffians. It has codimension three, and the Pfaffian multiplicity formula \cite[Theorem~2 and the example]{GK} gives
\[
             \deg Z=\frac{d(d+1)(2d+1)}{6}
                   =\frac{k(k+1)(k+2)}{24}.
\]
If every invertible input attained rank \(k\), Proposition~\ref{prop:degree} would give \(\gcd(k,\binom{k}{2},\binom{k}{3})\mid\deg Z\). Write \(a=v_2(k)\ge2\), where \(v_2\) is the exponent of two in an integer. The three binomial coefficients have respective valuations \(a,a-1,a\), since \(k-1\) is odd and \(v_2(k-2)=1\). Their gcd thus has valuation \(a-1\). In the degree formula \(k+1\) is odd and \(v_2(k+2)=1\), giving \(v_2(\deg Z)=a-2\). This contradicts divisibility.
\end{proof}

For \(k=4\), the locus \(Z\) is the Pl\"ucker Grassmannian \(\Gr(2,5)\subset\mathbb P^9\), and divisibility would require \(2\mid5\). At \(k=8\), it would require \(4\mid30\).

\subsection{Primitive cores of complex rank-four preservers}
Over \(\C\), Huang and Landsberg classify the four basic spaces of upper rank four \cite[Theorem~3.1]{HL}. Basic spaces satisfy the additional indecomposability and maximality conditions introduced by Eisenbud and Harris \cite{EH}. By \cite[Proposition~2.5]{HL}, every primitive space is obtained from basic spaces by rank-preserving parameter restrictions and output projections or by taking direct sums with ranks adding, followed by such operations.

To pass from spaces to maps, suppose that a model \(S\) of upper rank \(r\) gives a core by restricting to a subspace \(S'\subseteq S\) and projecting \(Y\mapsto LYR\). Any coefficient map into the resulting core lifts to \(S'\) by a linear section of this projection onto its image. At each input on which the core has rank \(r\), the lifted matrix must already have rank \(r\). Repeating this argument lifts through a sequence of projections. For a direct sum of total upper rank \(r=\sum_i r_i\), each lifted summand must attain its upper rank \(r_i\) on the same input. Thus it suffices to force a rank drop in one parent model or one summand.

For the two six-parameter models, put \(J(s,s')=\left(\begin{smallmatrix}s&s'\\0&s\end{smallmatrix}\right)\). For \(U,V,W\in\M_2\), set
\[
 H(U,V,W)=\begin{pmatrix}U&0&-W\\0&U&V\\V&W&0\end{pmatrix}.
\]
The two models are then
\begin{equation}\label{eq:hmodels}
\begin{aligned}
 \mathcal H_3&=\{H(sI_2,tI_2,D):s,t\in\C,\ D\in\M_2(\C)\},\\
 \mathcal H_4&=\{H(J(s,s'),J(t,t'),J(w,w')):s,s',t,t',w,w'\in\C\}.
\end{aligned}
\end{equation}

\begin{theorem}\label{thm:four}
Let \(\F=\C,\ k=4,\ t=\max(m,n)\). A bounded-image preserver in \eqref{eq:AL} has one of the following possible cores.
\begin{enumerate}[label=\textup{(\roman*)}]
\item \(r=0\), giving a compression.
\item \(r=3,\ a+b=1\), and \(m=n=4\). Up to transposition the core has the form \(TW_3(z(X))\), where \(W_3(z)=z\wedge-\) on \(\C^4\), \(z\) is an invertible change of \(Xv\) or \(X\T u\), and \(\ker T\) contains no nonzero decomposable bivector. The extra row or column must satisfy \eqref{eq:three}.
\item \(r=4,\ a=b=0\). Up to equivalence, zero padding and transposition, the map is \eqref{eq:wedge} with \(k=4\). Its projection and coefficients are those of Theorem~\ref{thm:wedge} and Proposition~\ref{prop:coeff}, with \(\Gamma\) surjective. The source has size \(4\times4,4\times5\), or \(5\times4\).
\end{enumerate}
Together with \eqref{eq:three} in (i) and (ii), these descriptions are sufficient. In particular, for \(t>5\) every complex rank-four preserver is standard with full factors or has compression form. The cutoff is attained by the full type II space on a \(4\times5\) source.
\end{theorem}
\begin{proof}
The core must attain upper rank \(r\) on every rank-four input by \eqref{eq:three}. A nonzero space of upper rank one cannot be primitive. A rank-two core \(\Alt_3\) has dimension three, contrary to Lemma~\ref{lem:codim}. A rank-three primitive core has at most four parameters by Atkinson's theorem, so that lemma gives \(m=n=4\). Its four parameters \(z\) must avoid zero on \(\GL_4\). Their common kernel is then a maximal singular space by Flanders's dimension bound, and its equality case gives the stated evaluations. Each evaluation attains every nonzero value on an invertible input. Hence \(T\) must be injective on \(z\wedge\C^4\) for every \(z\ne0\), as asserted in (ii).

Now take \(r=4\). A direct sum of lower positive primitive ranks must split as \(2+2\). The three coordinates of either \(\Alt_3\) summand vanish on some rank-four input by Lemma~\ref{lem:codim}. The lifting argument thus excludes direct sums. Table~\ref{tab:models} records the test for each basic model.

\begin{table}[hb]
\centering
\caption{The four basic complex spaces of upper rank four.}
\label{tab:models}
\smallskip
\small
\begin{tabularx}{\linewidth}{@{}l@{\enspace}l@{\enspace}c@{\enspace}c@{\quad}>{\raggedright\arraybackslash}X@{}}
\toprule
Type & Model & Size & \(\dim S\) & Test on rank-four inputs \\
\midrule
I & \(\Alt_5\) & \(5\times5\) & 10 &
  The degree of \(\Gr(2,5)\) is five and would have to be even. \\
II & \(z\wedge-\) & \(10\times5\) & 5 &
  Rank four persists under the two kernel conditions of Theorem~\ref{thm:wedge}. \\
III & \(\mathcal H_3\) & \(6\times6\) & 6 &
  A projective zero of a binary quadratic forces rank at most two. \\
IV & \(\mathcal H_4\) & \(6\times6\) & 6 &
  A subspace of codimension three consists of matrices of rank at most two. \\
\bottomrule
\end{tabularx}
\end{table}

\emph{Type I.} On a coordinate \(\M_4\) in the source, Corollary~\ref{cor:pfaffian} gives an invertible input where the parent alternating matrix has rank at most two.

\emph{Type III.} On the same coordinate block, write \(s,t,D\) for the six linear coordinates of the lifted map into \(\mathcal H_3\). When \(s=t=0\) its rank is \(2\rk D\). For each row \(D_i\) of \(D\), vanishing of the four entries of \((s,t,D_i)\) gives \(\rk D\le1\) and output rank at most two. These four entries cannot vanish together on \(\GL_4\). Their common kernel is a space of singular matrices of dimension at least 12. Flanders's theorem gives equality and shows that its annihilator consists of the forms \(X\mapsto u\T Xv\) with either \(u\) or \(v\) fixed. Two distinct four-dimensional spaces of this kind meet in at most one dimension. The two annihilators here contain the independent forms \(s,t\), since each quadruple has rank four. Thus they coincide, and all six coordinates factor through one evaluation \(X\mapsto Xv\) or \(X\mapsto X\T u\). Every nonzero value of the evaluation occurs on \(\GL_4\). In the target of this evaluation, the equations \(s=t=0\) define a two-plane. Over \(\C\), the quadratic form \(\det D\) vanishes at some nonzero point of this plane. The lifted matrix then has rank at most two.

\emph{Type IV.} At \(s=t=w=0\) all three blocks of \(H\) are multiples of \(J(0,1)\). The matrix is the Kronecker product of \(J(0,1)\) with a \(3\times3\) matrix of determinant zero, up to permuting factors. Its rank is at most two. The equations \(s=t=w=0\) define a subspace of codimension three, whose inverse image contains a rank-four input by Lemma~\ref{lem:codim}.

\emph{Type II.} The five parameters \(z\) must avoid zero on all rank-four inputs. If a projection restricted the five-dimensional domain of \(z\wedge-\) to a proper subspace contained in a hyperplane \(H\), the rank would drop whenever \(z\in H\). Lemma~\ref{lem:codim} gives a rank-four input in the inverse image of \(H\). Hence only projections on the exterior-square side can occur. Theorem~\ref{thm:wedge} and Proposition~\ref{prop:coeff} characterize these projections and coefficient maps. Proposition~\ref{prop:wprimitive} gives the surjectivity condition for primitivity.

The block criterion gives sufficiency in (i) and (ii), and Theorem~\ref{thm:wedge} gives it in (iii). On a \(4\times5\) source, the construction following
Proposition~\ref{prop:coeff} has the full primitive type II space as its image.
This space cannot lie in a compression of upper rank four, since a compression row or column could then be deleted to reduce the upper rank. Some rank-one input has rank-four output, which also excludes the standard forms.
\end{proof}

\section{Applications and open questions}\label{sec:apps}

\subsection{Rank distance and constant-rank codes}
With \(d(A,B)=\rk(A-B)\), a linear map is rank-\(k\) nonincreasing exactly when it takes pairs at distance at most \(k\) to pairs at distance at most \(k\). By Theorem~\ref{thm:main}, either it is nonexpansive at every distance or its image has diameter at most \(k\). Over a finite field a subspace of diameter at most \(k\) is a linear anticode. Thus Flanders's dimension bound gives \(\dim\im\phi\le k\max(p,q)\) in the second case. 
A standard map whose left factor has full column rank and whose right factor has full row rank preserves rank distance.
 Over \(\F_Q\), Proposition~\ref{prop:count} gives
\(Q^{mn}N_{n-i,q-j}\) ordered pairs
\((X,Y)\in\M_{m,n}(\F_Q)^2\) satisfying
\[
d(X,Y)=i,\qquad d\bigl(\phi(X),\phi(Y)\bigr)=j.
\]
Indeed, each difference \(Z=X-Y\) occurs for exactly \(Q^{mn}\)
ordered pairs, and \(\phi(X)-\phi(Y)=\phi(Z)\).
 Delsarte's bilinear-forms scheme and the code--anticode duality studied by Ravagnani give the coding-theoretic setting \cite{Delsarte,Ravagnani}. Maksaev, Medved and Promyslov treat a related two-map, bijective fixed-distance problem \cite{MMP}.

For square matrices over an arbitrary field, de Seguins Pazzis classified the linear maps \(\M_k\to\M_k\) that send \(\GL_k\) into itself \cite[Theorem~2]{Inv}. Besides the invertible standard maps, they are \(X\mapsto g(Xv)\) or \(g(X\T u)\), with \(u,v\ne0\) and a linear isomorphism \(g:\F^k\to W\subseteq\M_k\) onto a space whose nonzero members are invertible. For instance, take \(\F=\F_Q\) and let \(g(w)\) be multiplication by \(w\) on \(\F_{Q^k}\). The nonzero outputs have rank \(k\), giving a \(k\)-dimensional constant-rank code. Over an algebraically closed field, no such \(W\) has dimension \(k>1\), since the restriction of the determinant to \(W\) is a homogeneous form with a nontrivial zero.

For a source of arbitrary size, let \(v(X)\) be the column vector of the \(mn\) entries of \(X\). Form a \(kmn\times k\) block diagonal matrix with \(k\) diagonal blocks, each equal to \(v(X)\). Every nonzero input has output rank \(k\). This is a bounded-image preserver, and for \(k\ge2\) it increases the rank of each rank-one input. Li, Rodman and \v Semrl give this construction over \(\mathbb R\)
and \(\C\) \cite[Proposition~2.1]{LRS}. Their Proposition~2.2 gives
additional real examples using subspaces of nonsingular matrices.

\subsection{Bipartite tensors}
Identify \(\C^m\otimes\C^n\) with \(\M_{m,n}(\C)\) by \(u\otimes v\mapsto uv\T\). A vector has Schmidt rank \(j\) precisely when its coefficient matrix has rank \(j\). Suppose a complex linear operator between two such tensor products sends vectors of Schmidt rank at most \(k\) to vectors of Schmidt rank at most \(k\). If its image contains a vector of Schmidt rank greater than \(k\), Theorem~\ref{thm:main} shows that it acts by separate linear maps on the two tensor factors, possibly after interchanging them. Under our convention the two local factors of \(X\mapsto MXN\) are \(M\) and \(N\T\). Exact preservation in this case forces both local maps to be injective. This is the rectangular analogue of the linear preserver result used by Johnston in his study of separability measures \cite{Johnston}. It applies to state vectors and to an individual linear measurement operator before normalization.

For maps \(\Gamma\) with values in a single row of \(d\) entries, the least output dimension compatible with rank-one preservation depends on the field. Such a map preserves rank one exactly when the bilinear map \((u,v)\mapsto\Gamma(uv\T)\) is nonzero for every pair of nonzero vectors. The least attainable \(d\) is \(m+n-1\) over an algebraically closed field because the projective kernel must miss the Segre variety of dimension \(m+n-2\), and a general center of the required codimension does so. Over \(\F_Q\) the least \(d\) is \(\max(m,n)\), by Lemma~\ref{lem:codim} and multiplication of two linearly embedded subspaces of \(\F_{Q^{\max(m,n)}}\). This is the rank-one counterpart of Proposition~\ref{prop:wtarget}.

\subsection{Tests for exact preservation}
Over a finite field the joint distribution of the input and output ranks gives an exact test, even in cases such as Example~\ref{ex:binary}. Write \(\qbinom{a}{i}\) for the number of \(i\)-dimensional subspaces of \(\F_Q^a\), with value zero for \(i>a\).

\begin{proposition}\label{prop:count}
Let \(\F=\F_Q\) and \(k\le\min(p,q)\). For \(U\le\F^n,\ V\le\F^q\), compute
\[
 c(U,V)=\dim\{X:XU=0,\ \phi(X)V=0\},\qquad
 H_{ij}=\sum_{\dim U=i,\ \dim V=j} Q^{c(U,V)} .
\]
Define integers \(N_{ab}\), \(0\le a\le n,\ 0\le b\le q\), by the triangular system
\begin{equation}\label{eq:count}
 H_{ij}=\sum_{a\ge i,\ b\ge j}\qbinom{a}{i}\qbinom{b}{j} N_{ab}.
\end{equation}
Then \(\phi\) preserves rank \(k\) if and only if
\[
 N_{n-k,q-k}
 =\prod_{i=0}^{k-1}\frac{(Q^m-Q^i)(Q^n-Q^i)}{Q^k-Q^i}.
\]
\end{proposition}
\begin{proof}
Count the triples \((X,U,V)\) with \(U\subseteq\ker X\), \(V\subseteq\ker\phi(X)\) and prescribed subspace dimensions. Counting first over \(X\) gives \eqref{eq:count} with \(N_{ab}\) the number of inputs having the two nullities \(a,b\). Since the diagonal coefficients are one, solving in descending order of the two nullities gives every \(N_{ab}\). The product counts rank-\(k\) inputs by full-rank factorizations \(X=AB\) modulo a change of basis in \(\F^k\). Equality holds precisely when all rank-\(k\) inputs have output rank \(k\). Each \(c(U,V)\) is the nullity of the homogeneous linear system defined by \(XU=0\) and \(\phi(X)V=0\).
\end{proof}

Over an algebraically closed field, a bounded-image map can also be tested by ideal membership. Let \(I_j(Z)\) denote the ideal generated by the \(j\)-minors of \(Z\) in the polynomial ring in the source coordinates, with the zero ideal when there are no such minors. Put
\[
              J=I_{k+1}(X)+I_k(\phi(X)).
\]
Its common zeros are the inputs of rank at most \(k\) whose outputs have rank less than \(k\). Hence for \(\urk(\im\phi)\le k\), exact preservation is equivalent to \(I_k(X)\subseteq\sqrt J\). Equivalently, for every source \(k\)-minor \(\Delta\), the ideal \((J,1-z\Delta)\) in the polynomial ring with one extra variable \(z\) is the unit ideal. This follows from the Nullstellensatz and gives a finite Gr\"obner-basis test in concrete examples \cite{BV,Hartshorne}. For projected exterior products, the choice of field affects the test. Proposition~\ref{prop:wtarget} shows that a projection can preserve rank over a finite field and fail to preserve it over the algebraic closure.

\subsection{Open questions}
 
\begin{enumerate}
\item \emph{Exceptional rank-\(k\) preservers over \(\F_2\).}
Fix \(k\ge2\) and \(p,q\ge k+1\). For which source dimensions
\(m,n\ge k\) does there exist a linear rank-\(k\) preserver
\(\phi:\M_{m,n}(\F_2)\to\M_{p,q}(\F_2)\) such that
\(\rk\phi(X)>k\) for some \(X\) of rank less than \(k\)?
Example~\ref{ex:binary} gives such a map for \(m=n=k\) and every
such target by adjoining zero rows and columns to the output.
Classify these maps up to source and target equivalence, beginning
with the minimal case \(\M_k(\F_2)\to\M_{k+1}(\F_2)\).

\item 
\emph{The Pfaffian obstruction.} Does there exist a complex linear map \(\Gamma:M_6(\mathbb C)\to\operatorname{Alt}_7(\mathbb C)\) such that \(\operatorname{rank}\Gamma(X)=6\) for every \(X\in\mathrm{GL}_6(\mathbb C)\)? For \(k=6\), the divisor in the degree obstruction is \(\gcd(6,15,20)=1\), so the obstruction gives no restriction. More generally, for which \(k\ge6\) with \(k\equiv2\pmod4\) does such a map exist?
For a fixed candidate \(\Gamma\), let \(P_1,\ldots,P_7\) be the six-by-six principal Pfaffians of \(\Gamma(X)\). Then \(\Gamma\) has rank six on every invertible input if and only if \(\det X\in\sqrt{(P_1,\ldots,P_7)}\).

\end{enumerate}

\paragraph*{Use of generative AI.}
OpenAI's ChatGPT assisted with literature searches, revision of the
exposition, and the development of some proofs and examples, including the
proof of Lemma~\ref{lem:finitevanish} and the construction in
Example~\ref{ex:binary}.


\begin{thebibliography}{99}\small
\setlength{\itemsep}{.2\baselineskip}
\bibitem{Atkinson} M. D. Atkinson, Primitive spaces of matrices of bounded rank. II, \emph{J. Austral. Math. Soc. Ser. A} \textbf{34} (1983), 306--315. \doiurl{10.1017/S1446788700023740}.
\bibitem{ALLarge} M. D. Atkinson and S. Lloyd, Large spaces of matrices of bounded rank, \emph{Quart. J. Math. Oxford Ser. (2)} \textbf{31} (1980), 253--262. \doiurl{10.1093/qmath/31.3.253}.
\bibitem{AL} M. D. Atkinson and S. Lloyd, Primitive spaces of matrices of bounded rank, \emph{J. Austral. Math. Soc. Ser. A} \textbf{30} (1981), 473--482. \doiurl{10.1017/S144678870001795X}.
\bibitem{Beasley} L. B. Beasley, Linear operators on matrices: the invariance of rank-\(k\) matrices, \emph{Linear Algebra Appl.} \textbf{107} (1988), 161--167. \doiurl{10.1016/0024-3795(88)90242-X}.
\bibitem{BL} L. B. Beasley and T. J. Laffey, Linear operators on matrices: the invariance of rank-\(k\) matrices, \emph{Linear Algebra Appl.} \textbf{133} (1990), 175--184. \doiurl{10.1016/0024-3795(90)90248-B}.
\bibitem{Botta} P. Botta, Linear maps preserving rank less than or equal to one, \emph{Linear Multilinear Algebra} \textbf{20} (1987), 197--201. \doiurl{10.1080/03081088708817756}.
\bibitem{BV} W. Bruns and U. Vetter, \emph{Determinantal Rings}, Lecture Notes in Mathematics 1327, Springer, Berlin, 1988. \doiurl{10.1007/BFb0080378}.
\bibitem{Delsarte} P. Delsarte, Bilinear forms over a finite field, with applications to coding theory, \emph{J. Combin. Theory Ser. A} \textbf{25} (1978), 226--241. \doiurl{10.1016/0097-3165(78)90015-8}.
\bibitem{Inv} C. de Seguins Pazzis, The singular linear preservers of non-singular matrices, \emph{Linear Algebra Appl.} \textbf{433} (2010), 483--490. \doiurl{10.1016/j.laa.2010.03.021}.
\bibitem{Pazzis} C. de Seguins Pazzis, Local linear dependence seen through duality II, \emph{Linear Algebra Appl.} \textbf{462} (2014), 133--185. \doiurl{10.1016/j.laa.2014.08.013}.
\bibitem{Large} C. de Seguins Pazzis, Large spaces of bounded rank matrices revisited, \emph{Linear Algebra Appl.} \textbf{504} (2016), 124--189. \doiurl{10.1016/j.laa.2016.03.051}.
\bibitem{Binary} C. de Seguins Pazzis, Primitive spaces of matrices with upper rank two over the field with two elements, \emph{Linear Multilinear Algebra} \textbf{64} (2016), 1321--1353. \doiurl{10.1080/03081087.2015.1083526}.
\bibitem{DL} H. Duan and X. Lin, Topology of unitary groups and the prime orders of binomial coefficients, \emph{Sci. China Math.} \textbf{60} (2017), 1543--1548. \doiurl{10.1007/s11425-016-9004-0}.
\bibitem{EH} D. Eisenbud and J. Harris, Vector spaces of matrices of low rank, \emph{Adv. Math.} \textbf{70} (1988), 135--155. \doiurl{10.1016/0001-8708(88)90054-0}.
\bibitem{Flanders} H. Flanders, On spaces of linear transformations with bounded rank, \emph{J. London Math. Soc.} \textbf{37} (1962), 10--16. \doiurl{10.1112/jlms/s1-37.1.10}.
\bibitem{Fulton} W. Fulton, \emph{Intersection Theory}, 2nd ed., Ergebnisse der Mathematik und ihrer Grenzgebiete (3), vol.~2, Springer, Berlin, 1998. \doiurl{10.1007/978-1-4612-1700-8}.
\bibitem{GK} S. R. Ghorpade and C. Krattenthaler, The Hilbert series of Pfaffian rings, in \emph{Algebra, Arithmetic and Geometry with Applications}, Springer, Berlin, 2004, 337--356. \doiurl{10.1007/978-3-642-18487-1_22}.
\bibitem{GQ} R. Gow and R. Quinlan, On the vanishing of subspaces of alternating bilinear forms, \emph{Linear Multilinear Algebra} \textbf{54} (2006), 415--428. \doiurl{10.1080/03081080500361264}.
\bibitem{Hartshorne} R. Hartshorne, \emph{Algebraic Geometry}, Graduate Texts in Mathematics 52, Springer, New York, 1977. \doiurl{10.1007/978-1-4757-3849-0}.


\bibitem{HL} Hang Huang and J. M. Landsberg, On linear spaces of matrices of bounded rank, \emph{Selecta Math. (N.S.)} \textbf{32} (2026), article 30. \doiurl{10.1007/s00029-026-01137-x}.

\bibitem{HT} Huajun Huang and M.-C. Tsai, Linear rank nonincreasing and preserving maps between matrix spaces, 2026. \href{https://arxiv.org/abs/2610.40995}{arXiv:2610.40995}.


\bibitem{Johnston} N. Johnston, Characterizing operations preserving separability measures via linear preserver problems, \emph{Linear Multilinear Algebra} \textbf{59} (2011), 1171--1187. \doiurl{10.1080/03081087.2011.596540}.


\bibitem{LRS} C.-K. Li, L. Rodman and P. \v Semrl, Linear transformations between matrix spaces that map one rank specific set into another, \emph{Linear Algebra Appl.} \textbf{357} (2002), 197--208. \doiurl{10.1016/S0024-3795(02)00381-6}.
\bibitem{Loewy} R. Loewy, Linear mappings which are rank-\(k\) nonincreasing, \emph{Linear Multilinear Algebra} \textbf{34} (1993), 21--32. \doiurl{10.1080/03081089308818205}.
\bibitem{LoewyII} R. Loewy, Linear mappings which are rank-\(k\) nonincreasing, II, \emph{Linear Multilinear Algebra} \textbf{36} (1993), 115--123. \doiurl{10.1080/03081089308818281}.
\bibitem{MMP} A. Maksaev, N. Medved and V. Promyslov, Maps preserving a fixed rank-distance on matrices over finite fields, \emph{Finite Fields Appl.} \textbf{111} (2026), article 102759. \doiurl{10.1016/j.ffa.2025.102759}.
\bibitem{Ravagnani} A. Ravagnani, Rank-metric codes and their duality theory, \emph{Des. Codes Cryptogr.} \textbf{80} (2016), 197--216. \doiurl{10.1007/s10623-015-0077-3}.
\end{thebibliography}
\end{document}